\documentclass[a4paper]{amsart}

\usepackage{amsfonts}
\usepackage{amssymb}
\usepackage{amsmath}
\usepackage{hyperref}
\usepackage{mathrsfs}
\usepackage{centernot}
\usepackage{mathdots}
\usepackage{stmaryrd}
\usepackage[all]{xy}

\usepackage{mathtools}

\usepackage{color}

\theoremstyle{definition}

\theoremstyle{plain}

\newcommand{\rem}[1]{}

\newcommand{\N}{\mathbb{N}}

\newcommand{\Z}{\mathbb{Z}}

\newcommand{\frakp}{{\mathfrak{p}}}

\newcommand{\vphi}{\varphi}

\newcommand{\suchthat}{\,:\,}

\DeclareMathOperator{\Br}{Br} %
\DeclareMathOperator{\Cent}{Cent} %
\DeclareMathOperator{\End}{End} %
\DeclareMathOperator{\Hom}{Hom} %
\DeclareMathOperator{\id}{id} %
\DeclareMathOperator{\Nrd}{Nrd} %
\newcommand{\op}{\mathrm{op}} %
\DeclareMathOperator{\Pic}{Pic} %
\DeclareMathOperator{\Spec}{Spec} %
\DeclareMathOperator{\Sym}{Sym} %
\DeclareMathOperator{\Trd}{Trd} %

\newcommand{\nMat}[2]{\mathrm{M}_{#2}(#1)}

\newcommand{\trans}{{\mathrm{t}}}

\newcommand{\uO}{{\mathbf{O}}}

\newcommand{\uPGO}{{\mathbf{PGO}}}

\newcommand{\umu}{{\boldsymbol{\mu}}}

\newcommand{\fppf}{\mathrm{fppf}}

\newcommand{\units}[1]{{#1^\times}}

\usepackage[foot]{amsaddr}

\newtheorem{mainthm}{Theorem} 
\newtheorem{maincr}[mainthm]{Corollary} 
\newtheorem{mainlm}[mainthm]{Lemma}

\theoremstyle{definition}

\newtheorem{mainex}[mainthm]{Example} 
\newtheorem{mainrm}[mainthm]{Remark} 
\newtheorem*{ack}{Acknowledgement}

\title[Azumaya Algebras  With   Orthogonal Involution]{Azumaya Algebras  With   Orthogonal Involution Admitting an Improper Isometry}

\author{Uriya A.\ First$^*$}
\address{$^*$University of Haifa}
\email{uriya.first@gmail.com}

\begin{document}

\begin{abstract}
	Let $(A,\sigma)$ be an Azumaya algebra with   orthogonal involution over a ring $R$ with
	$2\in\units{R}$.
	We show that if $(A,\sigma)$ admits an improper
	isometry, i.e., an element $a\in A$ with $\sigma(a)a=1$ and $\Nrd_{A/R}(a)=-1$,
	then the Brauer class of $A$  is trivial. An   analogue of this statement
	also 
	holds for Azumaya algebras
	with quadratic pair  when $2\notin \units{R}$.
	We also show that at this level of generality, the  hypotheses
	do not guarantee that $A$ is a matrix algebra over $R$.
\end{abstract}

\maketitle

Let $R$ be a (commutative) ring 
and let $(A,\sigma,f)$ be an Azumaya algebra with a quadratic pair over
$R$. This means that $A$ is an Azumaya $R$-algebra of even degree,
$\sigma:A\to A$ is an orthogonal involution and $f$ is a semi-trace for $(A,\sigma)$,
i.e., an $R$-linear map from 
$\Sym(A,\sigma):=\{a\in A\suchthat \sigma(a)=a\}$ to $R$
satisfying $f(a+\sigma(a))=\Trd_{A/R}(a)$ for all $a\in A$.
See 
\cite[\S2.7, \S4.4]{Calmes_2015_groupes_classique} or \cite[\S4]{Gille_2023_Az_algs_quad_pairs_preprint} for all relevant definitions
and also \cite{Knus_1998_book_of_involutions} for the case where $R$ is a field.\footnote{
	Following \cite{Calmes_2015_groupes_classique} and \cite{Gille_2023_Az_algs_quad_pairs_preprint},
	we call an $R$-involution $\sigma$ on an Azumaya algebra $A$ orthogonal if there is a faithfully
	flat $R$-ring $R'$ splitting $A$ and such that $\sigma_{R'}:A_{R'}\to A_{R'}$
	is adjoint to a regular symmetric bilinear form over $R'$. We caution that this 
	is different from the definition used in \cite{Knus_1998_book_of_involutions}
	if $2\notin \units{R}$.
}
The isometry group of $(A,\sigma,f)$, denoted  $O(A,\sigma,f)$, is  the subgroup 
of $\units{A}$ consisting of elements
$a\in A$ satisfying $\sigma(a)a=1$ and $f(axa^{-1})=f(x)$ for all $x\in \Sym(A,\sigma)$.
The functor mapping an $R$-ring $S$ to the group $O(A_S,\sigma_S,f_S)$ (with $A_S=A\otimes_RS$, etc.) is
represented by an affine $R$-group scheme
$\uO(A,\sigma,f)$,
and there is a unique   $R$-group scheme homomorphism $\Delta$
from $\uO(A,\sigma,f)$ to the constant $R$-group scheme $(\Z/2\Z)_R$
whose kernel is (fiberwise over $\Spec R$) 
the identity connected component of $\uO(A,\sigma,f)$
\cite[4.4.0.43, 5.0.0.13, 2.7.0.32]{Calmes_2015_groupes_classique}.
As usual, elements of $O(A,\sigma,f)$ are called isometries, 
and an isometry $a$ is called \emph{proper}  if $\Delta(a)=0+2\Z$
and \emph{improper} if $\Delta(a)=1+2\Z$.

When $2\in\units{R}$, our setting  becomes simpler:
The map $f$ must coincide with
$\frac{1}{2}\Trd_{A/R}$, so we   omit it from the notation and
just write $O(A,\sigma)$ for the isometry group.
The elements of $O(A,\sigma)$ are the $a\in A$
satisfying $\sigma(a)a=1$, and if we identify
$(\Z/2\Z)_R$ with $\umu_{2,R}$ (the $R$-group scheme of square roots of $1$),
then the map $\Delta$ becomes the reduced norm $\Nrd_{A/R}:\uO(A,\sigma)\to \umu_{2,R}$.
An isometry $a\in O(A,\sigma)$ is therefore improper precisely when $\Nrd_{A/R}(a)=-1$.

Write $\Br R$ for the Brauer group of $R$ and let $[A]$
denote the Brauer class of $A$.
When $R$ is a field $F$ of characteristic not $2$, 
a result of Kneser \cite[Lem.~2.6.1b]{Kneser_1969_Galois_cohomology_classical_groups}
says that  $(A,\sigma,f)$ admits an improper
isometry if and only if $[A]=0$. 
A   proof working in any characteristic   was given later 
in \cite[Cor.~13.43]{Knus_1998_book_of_involutions}.
The statement was extended to the case where $R$  is a semilocal ring with $2\in\units{R}$ in \cite{First_2020_orthogonal_group},
where it was also shown that the ``if'' part of the statement is false for a general ring $R$.
Here we complete the picture by proving that the ``only if'' part holds for any ring $R$,
thus settling Question~3 in {\it op.\ cit.}

\begin{mainthm}\label{TH:main}
	Let $(A,\sigma,f)$ be an Azumaya algebra with a quadratic pair
	over a ring $R$. If $(A,\sigma,f)$ admits an improper isometry, then $[A]=0$ in $\Br R$.
\end{mainthm}

The proof is very different from the   argument in \cite{First_2020_orthogonal_group}.
It relies on the existence of generic Azumaya algebras with quadratic pair constructed implicitly
in \cite{First_2024_highly_versal_tors}, and earlier    
in \cite[\S5]{Auel_2019_Azumaya_algebras_without_inv} assuming $2\in \units{R}$,
in order to reduce the theorem to the known case
where $R$ is a field. 
It   does not generalize to Azumaya algebras with orthogonal involutions
over schemes, or more generally, locally ringed topoi, and so it remains open whether
Theorem~\ref{TH:main} holds in this broader context (consult \cite{Calmes_2015_groupes_classique}
and \cite[Dfn.~5.1, Ex.~7.4]{First_2020_involutions_of_Azumaya_algs}
for the   definitions  at this level of generality).

We   show in Example~\ref{EX:only-ex} below that
the conclusion   $[A]=0$
in Theorem~\ref{TH:main} cannot be improved to $A\cong \nMat{R}{n}$ for some $n\in\N$,
even if $R$ is connected.

\begin{mainrm}
	Theorem~\ref{TH:main} addresses only the case of even-rank Azumaya algebras.
	However, if $(A,\sigma)$ is an   Azumaya algebra with orthogonal involution of constant odd rank
	$n$,
	then $2[A]=0$ because $A\cong A^\op$ and $n[A]=0$ by a theorem of Saltman \cite{Saltman_1981_Brauer_group_is_torsion}, so we must have $[A]=0$.
\end{mainrm}

\begin{ack} 
An earlier version of this work  
treated only rings in which $2$ is invertible. 
We are grateful to an anonymous referee for insisting that we prove
Theorem~\ref{TH:main} for all rings, and giving many helpful suggestions toward achieving this.
\end{ack}

\section*{Proof of Theorem~\ref{TH:main}}

In what follows, an $R$-ring means a commutative $R$-algebra 
and
an $R$-domain is an $R$-ring which is moreover an integral domain.

We first prove the following result, which may be of independent interest.
We use the same conventions about torsors as those in \cite[\S2]{First_2024_highly_versal_tors}.

\begin{mainthm}\label{TH:smooth}
	Let $R_0$  be a noetherian   ring
	and let $G$ be a  linear $R_0$-group scheme    which is an extension
	of a    finite locally free $R_0$-group scheme by a reductive $R_0$-group scheme.
	Then every $G$-torsor over an affine $R_0$-scheme 
	arises as the base-change of a $G$-torsor
	over a smooth   $R_0$-scheme with connected geometric fibers.
\end{mainthm}

\begin{proof}
	Let $R$ be an $R_0$-ring and let $E$ be a $G$-torsor over $  R$.
	We need to show that there is a smooth  
	$R_0$-ring $R'$, a $G$-torsor $E'$ over $R'$
	and a morphism $f:\Spec R\to \Spec R'$ such that $E\cong E'\times_{\Spec R'}\Spec R$
	as $G$-torsors over $R$ and   $\Spec R'\to \Spec R_0$ has   connected geometric fibers.
	Write $R$ as a direct limit of its finitely generated $R_0$-algebras $\{R_i\}_{i\in I}$.
	Since $G$-torsors over $R_i$ are classified by
	the first  \v{C}ech cohomology group $\mathrm{H}^1_{\fppf}(  R_i,G)$
	and since $\mathrm{H}^1_{\fppf}(-,G)$ commutes with direct limits of rings
	\cite{Margaux_2007_limit_non_abel_coh}
	(see also \cite[Cor.~5.9, Rem.~5.14a]{SGA4_vol2}),
	there is $i\in I$ and a $G$-torsor $E_i$ over $R_i$
	such that $E $ is the base change of $E_i$ along $R_i\to R$.
	Replacing $R$ with $R_i$ and $E$ with $E_i$, we may assume   that $R$ is finitely generated over $R_0$,
	hence noetherian. In particular, $R$ has finite  Krull dimension.
	The existence of $E'$ and $R'$
	now follows from \cite[Thm.~8.1]{First_2024_highly_versal_tors}.
\end{proof}

\begin{maincr}\label{CR:smooth}
	Let $(A,\sigma,f)$ be an Azumaya algebra with quadratic pair of constant
	degree $2n$ over a ring
	$R$. 
	Then there is a   smooth   $\Z$-domain $S$,
	an Azumaya $S$-algebra with quadratic pair $(B,\tau,g)$ and a homomorphism
	$\vphi:R\to S$ such that $(B_R ,\tau_R,g_R)\cong (A,\sigma,f)  $.
\end{maincr}

\begin{proof} 
	Let $\uPGO_{2n}$ denote the automorphism $\Z$-group scheme of the split   Azumaya algebra
	with quadratic pair $(\nMat{\Z}{2n},\eta_{2n},f_{2n})$ over $\Z$ constructed in \cite[p.~57]{Calmes_2015_groupes_classique} or \cite[Ex.~4.5b]{Gille_2023_Az_algs_quad_pairs_preprint}; it is an extension of
	the finite contant $\Z$-group scheme $(\Z/2\Z)_\Z$
	by a reductive $\Z$-group scheme \cite[4.4.0.37, 8.1.0.55]{Calmes_2015_groupes_classique}.
	By \cite[4.4.0.34]{Calmes_2015_groupes_classique}, there is an equivalence
	of fibered categories over $\Spec\Z$ between the   $\uPGO_{2n}$-torsors
	and the  degree-$2n$ Azumaya algebras with quadratic pair.
	With this at hand,
	the corollary is just Theorem~\ref{TH:smooth}
	in the special case $R_0=\Z$ and $G=\uPGO_{2n}$. The resulting $R_0$-ring $S$
	is a domain by the following lemma.
\end{proof}

\begin{mainlm}\label{LM:connected}
	Let $f:X\to Y$ be a smooth morphism of schemes.
	If $Y$ is irreducible (resp.\ integral) and $f$ has connected fibers, then $X$ is also 
	irreducible
	(resp.\ integral).
\end{mainlm}

\begin{proof}
	The morphism $f$ is smooth, hence open
	\cite[\href{https://stacks.math.columbia.edu/tag/056G}{Tag 056G}]{DeJong_2020_stacks_project}.
	For $y\in Y$, let $X_y$ denote the scheme-theoretic fiber of $f$ over $y$,
	i.e., $X\times_{ Y}\Spec \kappa(y)$, where $\kappa(y)$ is the residue field of $y$;
	we use similar notation for open subschemes of $X$.
	Our assumptions on $f$ imply that $X_y$ is a smooth connected 
	$\kappa(y)$-scheme, hence integral.
	Note also that  $X_y$ is subspace of $X$ when both are viewed as topological spaces.
	
	We first prove that $X$ is irreducible.
	Let $U,U'$ be two nonempty open subsets of $X$. Then $f(U)$, $f(U')$
	are nonempty open subsets of the irreducible scheme $Y$, hence $f(U)\cap f(U')\neq \emptyset$.
	Let $y\in f(U)\cap f(U')$. Then $U_y$ and $U'_y$ are nonempty open subschemes of $X_y$.
	Since $X_y$ is irreducible, $U_y\cap U'_y\neq \emptyset$, so $U\cap U\neq \emptyset$.
	This proves that $X$ is irreducible.
	
	Suppose now that $Y$ is integral and let $y$ be its generic point.  Choose an open affine covering $\{V_i\}_i$ of $Y$
	and, for each $i\in I$, an open affine covering $\{U_{ij}\}_j$
	of $f^{-1}(V_i)$. Then each $U_{ij}\to V_i$ is smooth, and so corresponds
	to a smooth ring map $B_i\to A_{ij}$, where $B_i$ is a domain with fraction field
	$\kappa(y)$.
	Since $B_i\to A_{ij}$ is flat, $A_{ij}\to \kappa(y)\otimes_{B_i} A_{ij}$ is injective.
	Observe that $\Spec (\kappa(y)\otimes_{B_i}A_{ij})$ is an open subscheme of $X_y$,
	which is integral, so $\kappa(y)\otimes_{B_i} A_{ij}$ is a domain.
	This means that each $A_{ij}$ is a domain, and in particular reduced.
	We have therefore shown that $X$ is covered by reduced open affine subschemes, so $X$
	is reduced. Since $X$ is also irreducible,  it is integral.
\end{proof}

\begin{mainrm}
	When $2\in\units{R}$, we can give a shorter, more direct proof
	of Corollary~\ref{CR:smooth} using the generic Azumaya algebras with involution
	of \cite{Auel_2019_Azumaya_algebras_without_inv} as follows: We may forget about 
	the semi-trace $f$.
	Regard  $R$ as an algebra over $\Omega=\Z[\frac{1}{2}]$. 
	Let $m\in \N$ denote the \emph{Formanek number} of $A$ as defined in  
	\cite[\S5.1]{Auel_2019_Azumaya_algebras_without_inv}, let $(B,\tau)$ denote the 
	Azumaya algebra   with   orthogonal involution
	$(A(m,2n),\trans)$ constructed in {\it op.\ cit.} using the coefficient
	ring $\Omega$, and let $S=\Cent(B)$.
	By \cite[Thm.~17]{Auel_2019_Azumaya_algebras_without_inv} (see also
	\cite[Cor.~2.9b]{Saltman_1999_lectures_on_div_alg}), $S$ is a domain
	and there
	is a ring homomorphism $\phi:S\to R$ such that $(A,\sigma)\cong (B_R,\tau_R)$.
	Furthermore, by \cite[Prop.~20]{Auel_2019_Azumaya_algebras_without_inv},
	$S$ is a smooth $\Omega$-algebra, hence smooth over $\Z$.
\end{mainrm}

\begin{mainlm}\label{LM:spec}
	Let $R$, $(A,\sigma,f)$ be as in Corollary~\ref{CR:smooth}
	and let $a $ be an improper isometry of $(A,\sigma,f)$.
	Then there is an  smooth   $\Z$-domain $S$,
	an Azumaya $S$-algebra with quadratic pair $(B,\tau,g)$ admitting an improper
	isometry $b $, and a morphism
	$\vphi:R\to S$ such that $(B_R ,\tau_R,g_R)\cong (A,\sigma,f) $
	and the isomorphism maps $b\otimes 1$ to $a$. 
\end{mainlm}

\begin{proof}
	We  apply Corollary~\ref{CR:smooth} to
	$(A,\sigma,f)$
	to get a smooth  $\Z$-domain $S$ and an Azumaya $S$-algebra
	with quadratic pair $(B,\tau,g)$ such that
	$(A,\sigma,f)\cong (B_R,\tau_R,g_R)$ (it will not be our final $S$).
	
	Let $\uO^+(B, \tau,g)$ denote the $S$-group scheme of proper isomerties
	of $(B, \tau,g)$; it is semisimple \cite[8.1.0.55]{Calmes_2015_groupes_classique} and therefore smooth and
	has connected geometric fibers over $S$.
	We further let $\uO^-(B,\tau,g)$
	be the pullback of $\Delta:\uO^+(B, \tau,g)\to (\Z/2\Z)_S$
	along the $(1+2\Z)$-section $u:\Spec S\to (\Z/2\Z)_S$. Then $\uO^-(B,\tau,g)$
	the   affine $S$-group scheme representing the functor mapping an $S$-ring
	$T$ to the set of improper isometries of
	$(B_T, \tau_T,g_T)$, denoted $O^-(B_T,\tau_T,g_T)$. 
	The product in the group $O(B_T,\tau_T,g_T)$
	restricts to an action of $O^+(B_T,\tau_T,g_T)$
	on $O^-(B_T,\tau_T,g_T)$, which is free provided
	$O^-(B_T,\tau_T,g_T)\neq \emptyset$.
	By \cite[4.4.0.46, 4.4.0.37]{Calmes_2015_groupes_classique},
	there is a faithfully flat \'etale $S$-algebra
	$T$ with $O^-(B_T,\tau_T,g_T)\neq \emptyset$,
	so  $\uO^-(B,\tau,g)$ is a $\uO^+(B, \tau,g)$-torsor over $S$.
	Since $\uO^+(B,\tau,S)\to\Spec S$
	is smooth and has connected geometric fibers, the
	same holds for $\uO^-(B,\tau,g)\to \Spec S$. Our choice of $S$ and   Lemma~\ref{LM:connected}
	now imply that $\uO^-(B,\tau,g)$ is an integral scheme that is smooth over $\Spec \Z$.
	
	Let $S'$ denote the coordinate ring of the affine scheme $\uO^-(B,\tau,g)$;
	it is a smooth $\Z$-domain by what we have shown. 
	Put $B'=B_{S'}$, $\tau'=\tau_{S'}$, $g'=g_{S'}$. The identity map 
	$ S'\to   S'$ corresponds to a   improper isometry $b'\in O^-(B',\tau',g')$, which may
	be thought of as the universal improper isometry of $(B,\tau,g)$.
	Since $(A,\sigma,f)\cong (B_R,\tau_R,g_R)$, the improper isometry
	$a\in O^-(A,\sigma,f)=\uO^-(B,\tau,g)(R)$ corresponds to  an $S$-ring homomorphism 
	$\alpha :S'\to R$. Viewing $R$ as an $S'$-algebra
	via $\alpha $, we get $ (B',\tau',g')\otimes_{S'}R\cong 
	(B,\tau,g)\otimes_S R\cong
	(A,\sigma,f)$. Moreover, the equality $\alpha\circ \id_{S'}=\alpha$
	translates into the fact the image of $b\in O^-(B,\tau',g')$
	under the map $O^-(B,\tau',g')\to O^-(A,\sigma,f)$ is $a$.
	Thus, the data $\alpha:S'\to R$, $(B',\tau',g')$, $b'$  are what we were looking for.
\end{proof}

We can now deduce  Theorem~\ref{TH:main} from its known special case  
where $R$ is a field. 

\begin{proof}[Proof of Theorem~\ref{TH:main}]
	Recall that $\deg A$ denotes   the function from $\Spec R$ to $\N$
	mapping $\frakp$ to $\sqrt{\dim_{\kappa(\frakp)} A_{\kappa(\frakp)}}$,
	where $\kappa(\frakp)$ is the fraction field of $R/\frakp$.
	If $\deg A$ is not constant, then we can write $R$ as a product of rings $R=\prod_{i=1}^t R_i$
	such that $\deg A_{R_i}$ is constant for all $i$. It is enough to prove that $[A_{R_i}]=0$
	in $\Br R_i$ for each $i$, so we may work with each factor separately. We thus restrict
	to the case where
	$ \deg A$ is constant.

	Let $a$ be an improper isometry of $(A,\sigma,f)$
	and let $\vphi:S\to R$, $(B,\tau,g)$ and $b$ be as in Lemma~\ref{LM:spec}.
	Since $\vphi_*:\Br S\to \Br R$ maps $[B]$ to $[A]$, it is enough to prove
	that $[B]=0$.
	
	By construction, $S$ is a regular domain.
	Let $K$ be the fraction field of $S$.
	By the Auslander--Goldmann--Grothendieck Theorem \cite[Thm.~7.2]{Auslander_1960_Brauer_Group}
	(here we need $S$ to be regular), the map
	$\Br S\to \Br K$ is injective,
	so it   enough to show that
	$[B_K]=0$ in $\Br K$.
	But this follows from  \cite[Lem.~2.6.1b]{Kneser_1969_Galois_cohomology_classical_groups}
	and
	\cite[Cor.~13.43]{Knus_1998_book_of_involutions},
	because  $(B_K,\tau_K,g_K)$ is a central simple $K$-algebra with quadratic pair   admitting
	an improper isometry. 
\end{proof}

\begin{mainex}\label{EX:only-ex}
	We now give an example of a connected ring $R$ with $2\in\units{R}$
	and an Azumaya $R$-algebra with an orthogonal involution $(A,\sigma)$
	having an improper isometry, but such that $A\ncong \nMat{R}{n}$
	for all $n\in\N$.
	
	Let $R$ be a Dedekind domain with $2\in \units{R}$ and $\Pic R\cong \Z/2\Z$,
	and let $L$ be a  rank-$1$ projective $R$-module representing the nontrivial element of $\Pic R$;
	such rings $R$ exist, e.g., use \cite{Eakin_1973_class_groups}.
	In fact, any integral domain $R$ admitting a  rank-$1$ projective module
	$L$  such that 
	$2[L]=0$ in $\Pic R$ and $[L]\notin 2\Pic R$ will work.
	The former condition implies that there exists an $R$-module
	isomorphism $\phi:L\otimes_R L\to R$. 
	
	Put $M=R\oplus L$ (we write elements of $M$
	as column vectors), $A=\End_R(M)$, and define
	a symmetric $R$-bilinear form $b:M\times M\to R$ by 
	$b([\begin{smallmatrix} x \\ y\end{smallmatrix}], [\begin{smallmatrix} z \\ w\end{smallmatrix}])
	=xz+\phi(y\otimes w)$. 
	If we were to choose $L= R$, then $b$ would be  a
	diagonal bilinear form $([\begin{smallmatrix} x \\ y\end{smallmatrix}], [\begin{smallmatrix} z \\ w\end{smallmatrix}])
	=xz+\alpha y w$, where $\alpha\in\units{R}$ depends on $\phi$.
	Since $L$ becomes isomorphic to $R$ over some Zariski covering of $R$,
	this means that the $R$-linear map  
	$\hat{b}:M\to M^*=\Hom_R(M,R)$ given by $(\hat{b}x)y=b(x,y)$ is an isomorphism
	Zariski-locally on $R$,
	hence an isomorphism;
	otherwise said, 
	$b$ is \emph{regular}. 
	Thus, $b$ is adjoint to an orthogonal involution $\sigma:A\to A$;
	for $a\in A$, the element $\sigma(a)$ is the unique $R$-endomorphism
	of $M$ satisfying   
	$b(ax,y)=b(x,\sigma(a) y)$
	for all $x,y\in M$.
	
	Let $a=\id_R\oplus (-\id_L)\in A$. Then $\sigma(a)=a$, hence $a\sigma(a)=a^2=1_A$
	and $a\in O(A,\sigma)$. Writing $F$ for the fraction field   of $R$,
	we have 
	\[\Nrd_{A/R}(a)=\Nrd_{\nMat{F}{2}/F}(\id_F\oplus (-\id_F))=-1,\] 
	so $a$ is an improper isometry of $(A,\sigma)$.

	We now show that $A\ncong \nMat{R}{n}$ as $R$-algebras for any $n\in\N$.
	For the sake of contradiction, suppose  that such an isomorphism exists; rank considerations
	then force $n=2$.
	Consider $R^2$ as a left $A$-module via the isomorphism
	$\End_R(R^2)\cong\nMat{R}{2}\cong A$
	and put $P=\Hom_A(R^2,M)$. We claim that $P$ is a projective $R$-module of rank $1$. That 
	$P$ is finitely generated and projective follows readily from the fact that 
	that $M$ and $R^2$ are projective left $A$-modules. 
	This also implies that $P\otimes_R k\cong \Hom_{A\otimes k}(k^2,M\otimes_R k)$
	for every $R$-field $k$.
	Since $A$ is an Azumaya $R$-algebra of degree $2$ with trivial Brauer class,
	$A\otimes_R k\cong \nMat{k}{2}$, and under this isomorphism $M\otimes_Rk\cong k^2$
	as left $\nMat{k}{2}$-modules. Thus, $\dim_k P\otimes_R k=\dim_k \End_{\nMat{k}{2}}(k^2)=1$,
	proving that $P$ is of rank $1$.
	Consider the morphism $\vphi:P\otimes_R R^2=\Hom_A(R^2,M)\otimes_R R^2\to M$ given by
	$\vphi(p\otimes x)=p(x)$ ($p\in P$, $x\in R^2$). It is an isomorphism
	because the source and target are finitely generated projective $R$-modules,
	and because $\vphi\otimes_R \id_k$ is an isomorphism for any $R$-field $k$ by our earlier observations.
	Thus, $P^2\cong P\otimes_R R^2\cong M=R\oplus L$ as $R$-modules.
	Taking the second exterior power of both sides, we find that $P\otimes_R P\cong L$,
	or rather, $2[P]=[L]$ in $\Pic R$. This contradicts our choice of $L$, so an
	isomorphism $A\to \nMat{R}{n}$ cannot exist.
\end{mainex}

\bibliographystyle{plain}
\bibliography{MyBib_24_03}

\def\cprime{$'$} \def\cprime{$'$} \def\cprime{$'$} \def\cprime{$'$}
\begin{thebibliography}{10}

\bibitem{SGA4_vol2}
{\em Th\'eorie des topos et cohomologie \'etale des sch\'emas. {T}ome 2},
  volume Vol. 270 of {\em Lecture Notes in Mathematics}.
\newblock Springer-Verlag, Berlin-New York, 1972.
\newblock S\'eminaire de G\'eom\'etrie Alg\'ebrique du Bois-Marie 1963--1964
  (SGA 4), Dirig\'e{} par M. Artin, A. Grothendieck et J. L. Verdier. Avec la
  collaboration de N. Bourbaki, P. Deligne et B. Saint-Donat.

\bibitem{Auel_2019_Azumaya_algebras_without_inv}
Asher Auel, Uriya~A. First, and Ben Williams.
\newblock Azumaya algebras without involution.
\newblock {\em J. Eur. Math. Soc. (JEMS)}, 21(3):897--921, 2019.

\bibitem{Auslander_1960_Brauer_Group}
Maurice Auslander and Oscar Goldman.
\newblock The {B}rauer group of a commutative ring.
\newblock {\em Trans. Amer. Math. Soc.}, 97:367--409, 1960.

\bibitem{Calmes_2015_groupes_classique}
Baptiste Calm\`es and Jean Fasel.
\newblock Groupes classiques.
\newblock In {\em Autours des sch\'emas en groupes. {V}ol. {II}}, volume~46 of
  {\em Panor. Synth\`eses}, pages 1--133. Soc. Math. France, Paris, 2015.

\bibitem{Eakin_1973_class_groups}
Paul Eakin and W.~Heinzer.
\newblock More noneuclidian {${\rm PID}$}'s and {D}edekind domains with
  prescribed class group.
\newblock {\em Proc. Amer. Math. Soc.}, 40:66--68, 1973.

\bibitem{First_2020_orthogonal_group}
Uriya~A. First.
\newblock On the non-neutral component of outer forms of the orthogonal group.
\newblock {\em J. Pure Appl. Algebra}, 225(1):106477, 6, 2021.

\bibitem{First_2024_highly_versal_tors}
Uriya~A. First.
\newblock Highly versal torsors.
\newblock In {\em Amitsur {C}entennial {S}ymposium}, volume 800 of {\em
  Contemp. Math.}, pages 129--174. Amer. Math. Soc., [Providence], RI, [2024]
  \copyright2024.

\bibitem{First_2020_involutions_of_Azumaya_algs}
Uriya~A. First and Ben Williams.
\newblock Involutions of {A}zumaya algebras.
\newblock {\em Doc. Math.}, 25:527--633, 2020.

\bibitem{Gille_2023_Az_algs_quad_pairs_preprint}
Philippe Gille, Erhard Neher, and Cameron Ruether.
\newblock Azumaya algebras and obstructions to quadratic pairs over a scheme.
\newblock 2023.
\newblock arXiv:2209.07107.

\bibitem{Kneser_1969_Galois_cohomology_classical_groups}
M.~Kneser.
\newblock {\em Lectures on {G}alois cohomology of classical groups}.
\newblock Tata Institute of Fundamental Research, Bombay, 1969.
\newblock With an appendix by T. A. Springer, Notes by P. Jothilingam, Tata
  Institute of Fundamental Research Lectures on Mathematics, No. 47.

\bibitem{Knus_1998_book_of_involutions}
Max-Albert Knus, Alexander Merkurjev, Markus Rost, and Jean-Pierre Tignol.
\newblock {\em The book of involutions}, volume~44 of {\em American
  Mathematical Society Colloquium Publications}.
\newblock American Mathematical Society, Providence, RI, 1998.
\newblock With a preface in French by J. Tits.

\bibitem{Margaux_2007_limit_non_abel_coh}
Benedictus Margaux.
\newblock Passage to the limit in non-abelian \v{C}ech cohomology.
\newblock {\em J. Lie Theory}, 17(3):591--596, 2007.

\bibitem{Saltman_1981_Brauer_group_is_torsion}
David~J. Saltman.
\newblock The {B}rauer group is torsion.
\newblock {\em Proc. Amer. Math. Soc.}, 81(3):385--387, 1981.

\bibitem{Saltman_1999_lectures_on_div_alg}
David~J. Saltman.
\newblock {\em Lectures on division algebras}, volume~94 of {\em CBMS Regional
  Conference Series in Mathematics}.
\newblock Published by American Mathematical Society, Providence, RI, 1999.

\bibitem{DeJong_2020_stacks_project}
The {Stacks Project Authors}.
\newblock \textit{Stacks Project}.
\newblock \url{https://stacks.math.columbia.edu}, 2020.

\end{thebibliography}

\end{document}